\documentclass[fleqn]{mat01}
\usepackage{times,mathtimy,amssymb,latexsym,amsbsy}
\begin{document}

\setcounter{page}{461}
\firstpage{461}

\newtheorem{theore}{Theorem}
\renewcommand\thetheore{\arabic{section}.\arabic{theore}}
\newtheorem{theor}[theore]{\bf Theorem}
\newtheorem{rem}[theore]{Remark}
\newtheorem{propo}[theore]{\rm PROPOSITION}
\newtheorem{lem}[theore]{Lemma}
\newtheorem{definit}[theore]{\rm DEFINITION}
\newtheorem{coro}[theore]{\rm COROLLARY}
\newtheorem{exampl}[theore]{Example}
\newtheorem{case}{Case}
\renewcommand\thecase{{\it \Roman{case}}}

\newtheorem{note}[theore]{Note}

\def\lemm{\trivlist \item[\hskip \labelsep{\it Lemma {\rm (}Fatou property{\rm )}.}]}

\def \rofi{\rho_{\phi}}
\def \cfi{ces_{\phi}}
\def \set#1{\{#1\}}


\font\chuto=cmbx12 at 12pt
\font\nho=cmr8
\font\dam=cmbx8
\font\vua=cmcsc10

\renewcommand\theequation{\thesection\arabic{equation}}

\title{Basic topological and geometric properties of Ces\`{a}ro--Orlicz spaces}

\markboth{Yunan Cui et~al}{Basic topological and geometric properties of Ces\`{a}ro--Orlicz spaces}

\author{YUNAN CUI$^{\rm a}$, HENRYK HUDZIK$^{\rm b}$, NARIN PETROT$^{\rm c}$, SUTHEP SUANTAI$^{\rm c}$ and ALICJA SZYMASZKIEWICZ$^{\rm d}$}

\address{$^{\rm a}$Department of Mathematics, Harbin University of Science and Technology, Harbin 150~080,
People's Republic of China\\
\noindent$^{\rm b}$Faculty of Mathematics and Computer Science,
Adam Mickiewicz University, Pozna\'{n}, Poland\\
\noindent$^{\rm c}$Department of Mathematics, Faculty of Science,
Chiang Mai University, Chiang Mai, Thailand\\
\noindent$^{\rm d}$Institut of Mathematics,
Technical University, Szczecin, Poland\\
\noindent E-mail: cuiya@mail.hrbust.edu.cn; hudzik@amu.edu.pl; npetrot@yahoo.com; scmti005@chiangmai.ac.th;
Alicja.Szymaszkiewicz@ps.pl}

\volume{115}

\mon{November}

\parts{4}

\pubyear{2005}

\Date{MS received 29 April 2005}

\begin{abstract}
Necessary and sufficient conditions under which the Ces\`{a}ro--Orlicz
sequence space $\cfi$ is nontrivial are presented. It is proved that for the
Luxemburg norm, Ces\`{a}ro--Orlicz spaces $\cfi$ have the Fatou property.
Consequently, the spaces are complete. It is also proved that the
subspace of order continuous elements in $\cfi$ can be defined in two
ways. Finally, criteria for strict monotonicity, uniform monotonicity
and rotundity (= strict convexity) of the spaces $\cfi$ are given.
\end{abstract}


\maketitle

\section{Introduction}

As usual, $\Bbb R$, $\Bbb R_+$ and $\Bbb N$ denote the sets of reals,
nonnegative reals and natural numbers, respectively. The space of all
real sequences $x=(x(i))_{i=1}^{\infty}$ is denoted by $l^0$.

A map \hbox{$\phi\!\!:\Bbb R\to [0,+\infty]$} is said to be an Orlicz function if
$\phi$ is even, convex, left continuous on $\Bbb R_+$, continuous at
zero, $\phi(0)=0$ and $\phi(u)\to\infty$ as $u\to\infty$. If $\phi$
takes value zero only at zero we will write $\phi>0$ and if $\phi$ takes
only finite values we will write $\phi<\infty$ \cite{1,13,17,18,19,20}.

The arithmetic mean map $\sigma$ is defined on $l^0$ by the formula:
\begin{equation*}
\sigma x=(\sigma x(i))_{i=1}^{\infty},\quad\text{where} \ \ \sigma x(i)=\frac 1 i\sum_{j=1}^{i}|x(j)|.
\end{equation*}
Given any Orlicz function $\phi$, we define on $l^0$ the following two
convex modulars \cite{18,19}
\begin{equation*}
I_{\phi}(x)=\sum_{i=1}^{\infty}\phi(x(i)),\quad \rofi(x)=I_{\phi}(\sigma x).
\end{equation*}
The space
\begin{equation*}
\cfi=\{x\in l^0:\quad \rofi(\lambda x)<\infty \ \ \text{for some} \ \ \lambda>0\},
\end{equation*}
where $\phi$ is an Orlicz function which is called the
Ces\`{a}ro--Orlicz sequence space. We equip this space with the Luxemburg norm
\begin{equation*}
\|x\|_{\phi}=\inf\left\{\lambda>0:\quad \rofi\left(\frac{x}{\lambda}\right)\leq 1\right\}.
\end{equation*}
In the case when $\phi(u)=|u|^p$, $1\leq p<\infty$, the space $\cfi$ is
nothing but the Ces\`{a}ro sequence space $ces_p$ (see \cite{5,6,7,14,16,21}) and
the Luxemburg norm generated by this power function is then expressed by
the formula
\begin{equation*}
\|x\|_{ces_p}=\left[\sum_{i=1}^\infty \left(\frac{1}{i}\sum_{j=1}^{i} |x(j)|\right)^{p}\right]^{\frac{1}{p}}.
\end{equation*}

\looseness -1 A Banach space $(X,\|\cdot\|)$ which is a subspace of $l^0$ is said to
be a K\"{o}the sequence space, if:

\begin{enumerate}
\renewcommand{\labelenumi}{(\roman{enumi})}
\item for any $x\in l^0$ and $y\in X$ such that $|x(i)|\leq|y(i)|$
for all $i\in \Bbb N$, we have $x\in X$ and $\|x\|\leq\|y\|$,

\item there is $x\in X$ with $x(i)\not=0$ for all $i\in\Bbb N$.
\end{enumerate}\vspace{-1pc}

Any nontrivial Ces\`{a}ro--Orlicz sequence space belongs to the class of K\"{o}the
sequence spaces.

An element $x$ from a K\"{o}the sequence space $(X,\|\cdot\|)$ is called
order continuous if for any sequence $(x_n)$ in $X_+$ (the positive cone
of $X$) such that $x_n\leq|x|$ and $x_n\to 0$ coordinatewise, we have
$\|x_n\|\to 0$.

A K\"{o}the sequence space $X$ is said to be order continuous if any
$x\in X$ is order continuous. It is easy to see that $X$ is order
continuous if and only if $\|(0,\dots,0,x(n+1),x(n+2),\dots)\|\to 0$ as
$n\to\infty$ for any $x\in X$.

A K\"{o}the sequence space $X$ is called monotone complete if for any
$x\in X_+$ and any sequence $(x_n)$ in $X_+$ such that $x_n(i)\leq
x_{n+1}(i)\leq\dots\leq x(i)$ for all $i\in\Bbb N$ and $x_n\to x$
coordinatewise, we have $\|x_n\|\to\|x\|$.

We say a K\"{o}the sequence space $X$ has the Fatou property if for any
sequence $(x_n)$ in $X_+$ and any $x\in l^0$ such that $x_n\to x$
coordinatewise and $\sup_ n\|x_n\|<\infty$, we have that $x\in X$ and
$\|x_n\|\to\|x\|$. For the above properties of K\"{o}the sequence (and
function) spaces we refer to [12] and [15].

We say an Orlicz function $\phi$ satisfies the $\Delta_2$-condition at
zero ($\phi\in\Delta_2(0)$ for short) if there are $K>0$ and $a>0$ such
that $\phi(a)>0$ and $\phi(2u)\leq K\phi(u)$ for all $u\in [0,a]$.

A modular $\rho$ (for its definition see \cite{4,18,19}) is said to
satisfy the $\Delta_2$-condition if for any $\epsilon>0$ there exist
constants $k\geq 2$ and $a>0$ such that $\rho(2x)\leq k\rho(x)+\epsilon$
for all $x\in X$ with $\rho(x)\leq a$.

If $\rho$ satisfies the $\Delta_2$-condition for any $a>0$ and
$\epsilon>0$ with $k\geq 2$ dependent on $a$ and $\epsilon$, we say that
$\rho$ satisfies the strong $\Delta_2$-condition ($\rho\in\Delta_2^{S}$
for short) (see [4]).

We say a K\"{o}the sequence space $X$ is strictly monotone, and then we
write $X\in (SM)$, if $\|x\|<\|y\|$ for all $x,y\in X$ such that $0\leq
x\leq y$ and $x\not =y$.

We say a K\"{o}the sequence space $X$ is uniformly monotone, and then we
write $X\in (UM)$, if for each $\epsilon>0$ there exists
$\delta(\epsilon)>0$ such that for any $x,y\geq 0$ such that $\|x\|=1$
and $\|y\|\geq\epsilon$, we have $\|x+y\|\geq1+\delta(\epsilon)$.

Let $B(X)$ (resp. $S(X)$) be the closed unit ball (resp. the unit
sphere) of $X$. A point $x\in S(X)$ is called an {\it extreme point} of
$B(X)$ if for every $y,z \in B(X)$ the equality $2x=y+z$ implies $y=z.$
Let Ext $B(X)$ denote the set of all extreme points of
$B(X)$. A Banach space $X$ is said to be {\it rotund} (write {\bf (R)}
for short), if Ext $B(X)= S(X)$. For these and other
geometric notions of rotundity type and their role in mathematics we
refer to the monographs \cite{1,8,19} and also to the papers \cite{2,3,10,11,22}.

We say that $u\in \Bbb R$ is a {\it point of strict convexity} of $\phi$
if $\phi\left(\frac{v+w}{2}\right) <\frac{\phi(v)+\phi(w)}{2},$ whenever
$u=\frac{v+w}{2}$ and $v\not = w$. We denote by $S_\phi$ the set of all
{\it points of strict convexity} of $\phi.$

An interval $[a, b]$ is called a {\it structurally affine interval} for
an Orlicz function $\phi$, or simply, SAI of $\phi$, provided that
$\phi$ is affine on $[a, b]$ for any $\varepsilon >0$ and it is not
affine either on $[a-\varepsilon, b]$ or on $[a, b+\varepsilon]$. Let
$\{[a_i, b_i]\}_i$ be all the SAIs of $\phi$. It is obvious that
\begin{equation*}
S_\phi = \Bbb R\, \big\backslash \bigcup_i \,(a_i, b_i).
\end{equation*}

\section{Results}

First we present necessary and sufficient conditions for nontriviality
of $\cfi$.

\begin{theor}[\!]
The following conditions are equivalent{\rm :}

\begin{enumerate}
\renewcommand{\labelenumi}{\rm \arabic{enumi})}
\item $\cfi\not =\{0\},$

\item $\exists_{n_1} \sum_{n=n_1}^{\infty}\phi\left(\frac{1}{n}\right)<\infty,$

\item $\forall_k>0\quad\exists_{n_k}\sum_{n=
n_k}^{\infty}\phi\left(\frac{k}{n}\right)<\infty.$
\end{enumerate}
\end{theor}

\begin{proof}$\left.\right.$\vspace{.5pc}

\noindent (1) $\Rightarrow$ (2). Let $0\not = z\in\cfi$. Since $z\not =0$, there
exists $l\in\Bbb N$ such that $z(l)\not =0$. Hence
$y=(0,\dots,0,z(l),0,\dots)\in\cfi$, and consequently,
$x=(0,\dots,0,1,0\dots)\in\cfi$, which means that there exists $k>0$ such
that $\rofi(kx)=\sum_{n=l}^{\infty}\phi\left(\frac{k}{n}\right)<\infty$.
We will consider two cases:

\begin{enumerate}
\renewcommand{\labelenumi}{\arabic{enumi}.}
\item $k>1$. Then for all $n$ we have $\frac 1 n<\frac k n$. From
monotonicity of the function $\phi$ we have $\phi(\frac{1}{n})<\phi(\frac{k}{n})$ for all $n$. Therefore
\begin{equation*}
\sum_{n=l}^{\infty}\phi\left(\frac{1}{n}\right)<\sum_{n=l}^{\infty}\phi\left(\frac{k}{n}\right)<\infty.
\end{equation*}
So it is enough to take $n_1=l$.

\item $0<k<1$. Then there exists $m\in \Bbb N$ such that $\frac 1
m\leq k$, whence $\frac 1 {m n}\leq\frac k n$ for all $n\in \Bbb N$ and
so, $\sum_{n=l}^{\infty}\phi\left(\frac{1}{mn}\right)\leq\sum_{n=l}^{\infty}\phi(\frac{k}{n})$. Consequently,
\begin{align*}
&\sum_{n=ml}^{\infty}\phi\left(\frac{1}{n}\right)=
\phi\left(\frac{1}{ml}\right)+\phi\left(\frac{1}{ml+1}\right)+\cdots+\phi\left(\frac{1}{ml+(m-1)}\right)\\[1pc]
&\quad\,+\phi\left(\frac{1}{m(l+1)}\right)+\phi\left(\frac{1}{m(l+1)+1}\right)\\[1pc]
&\quad\,+\cdots+\phi\left(\frac{1}{m(l+1)+(m-1)}\right)+\cdots\leq\phi\left(\frac{1}{ml}\right)+\phi\left(\frac{1}{ml}\right)
\end{align*}
\begin{align*}
&\quad\,+\cdots+\phi\left(\frac{1}{ml}\right)+\phi\left(\frac{1}{m(l+1)}\right)+\phi\left(\frac{1}{m(l+1)}\right)\\[.8pc]
&\quad\,+\cdots+\phi\left(\frac{1}{m(l+1)}\right)+\cdots=m\phi\left(\frac{1}{ml}\right)\\[.8pc]
&\quad\,+m\phi\left(\frac{1}{m(l+1)}\right)+\cdots=m\sum_{n=l}^{\infty}\phi\left(\frac{1}{mn}\right)\!\leq\! m\sum_{n=l}^{\infty}\phi\left(\frac{k}{n}\right)<\infty.
\end{align*}
Taking $n_1:=ml$, we get the thesis of condition (2).
\end{enumerate}

\noindent (2) $\Rightarrow$ (3). Assume that there exists $n_1$ such that
$\sum_{n=n_1}^{\infty}\phi\left(\frac{1}{n}\right)<\infty$ and consider two cases.

\begin{enumerate}
\renewcommand{\labelenumi}{\arabic{enumi}.}
\item $0<k<1$. Then $\frac k n<\frac 1 n$ and
$\sum_{n=n_1}^{\infty}\phi\left(\frac{k}{n}\right)<\sum_{n=n_1}^{\infty}\phi\left(\frac{1}{n}\right)
<\infty$. Taking $n_k:=n_1$, we have $\sum_{n=n_k}^{\infty}\phi\left(\frac{k}{n}\right)<\infty$.

\item $k>1$. Then there exists $m\in\Bbb N$ such that $k\leq m$.
Defining $n_k:=n_1 m$, we have\vspace{.5pc}
\begin{align*}
&\sum_{n=n_k}^{\infty}\phi\left(\frac{k}{n}\right)\leq\sum_{n=n_k}^{\infty}\phi\left(\frac{m}{n}\right)=\sum_{n=n_1m}^{\infty}\phi\left(\frac{m}{n}\right)=\phi\left(\frac{m}{n_1m}\right)\\[.5pc]
&\qquad\,+\phi\left(\frac{m}{n_1m+1}\right)+\cdots+\phi\left(\frac{m}{n_1m+(m-1)}\right)\!+\!\phi\left(\frac{m}{(n_1+1)m}\right)\\[.5pc]
&\qquad\,+\phi\left(\frac{m}{(n_1+1)m+1}\right)+\cdots+\phi\left(\frac{m}{(n_1+1)m+(m-1)}\right)+\cdots\\[.5pc]
&\quad\,\leq\phi\left(\frac{1}{n_1}\right)+\phi\left(\frac{1}{n_1}\right)+\cdots+\phi\left(\frac{1}{n_1}\right)\\[.5pc]
&\qquad\,+\phi\left(\frac{1}{n_1+1}\right)+\phi\left(\frac{1}{n_1+1}\right)+\cdots+\phi\left(\frac{1}{n_1+1}\right)+\cdots\\[.5pc]
&\quad\,=m\phi\left(\frac{1}{n_1}\right)+m\phi\left(\frac{1}{n_1+1}\right)+\cdots=m \sum_{n=n_1}^{\infty}\phi\left(\frac{1}{n}\right)<\infty.
\end{align*}
\end{enumerate}

\noindent (3) $\Rightarrow$ (1). Take $k=1$. By the assumption that condition (3)
holds, there exists $n_1\in \Bbb N$ such that
$\sum_{n=n_1}^{\infty}\phi\left(\frac{1}{n}\right)<\infty$. Define $x=(\underbrace
{0,\dots,0,}_{n_{_{1}}-1\text{ times}}1,0,\dots)$. Clearly, $x\in l^0$ and
\begin{equation*}
\rofi(kx)=\rofi(x)=\sum_{n=n_1}^{\infty}\phi\left(\frac{1}{n}\right)<\infty.
\end{equation*}
Hence $x\in \cfi$.$\hfill\Box$
\end{proof}

We will assume in the following that $\cfi$ is nontrivial, that is,
conditions (2) and (3) from Theorem 2.1 hold. Our next theorem gives some
sufficient conditions for the nontriviality of $\cfi$ in terms of some
lower index for the generating Orlicz function $\phi$.

\begin{theor}[\!]
For the conditions{\rm :}

\begin{enumerate}
\renewcommand{\labelenumi}{\rm (\alph{enumi})}
\item $\liminf_{t\to 0}\frac{t\phi'(t)}{\phi(t)}>1${\rm ,}

\item $\exists_{\epsilon>0}\exists_{A>0}\exists_{u_0>0} \ \forall_{0\leq u\leq u_0}\,\phi(u)\leq Au^{1+\epsilon}${\rm ,}

\item $\exists n_1
\sum_{n=n_1}^{\infty}\phi\left(\frac{1}{n}\right)<\infty${\rm
,}\vspace{-1pc}
\end{enumerate}
we have the implications ${\rm (a)}\Rightarrow {\rm (b)} \Rightarrow {\rm (c)}${\rm .}
\end{theor}

\begin{proof}$\left.\right.$\vspace{.4pc}

\noindent (a) $\Rightarrow$ (b). Although this implication appeared for example in
[9] we will present its proof for the sake of completeness.

By the assumption that $\liminf_{t\to 0}\frac{t\phi'(t)}{\phi(t)}>1$ we
know that there exists $t_0$ such that $\alpha:=\inf_{0<t\leq
t_0}\frac{t\phi'(t)}{\phi(t)}>1$. Then for all $0\leq t\leq t_0$ we have
that $\frac{t\phi'(t)}{\phi(t)}\geq \alpha$, that is,
$\frac{\phi'(t)}{\phi(t)}\geq \frac{\alpha}t$. Take $0<\lambda<1$. Then
$\lambda t<t$ and so for $0<t\leq t_0$:\vspace{-.3pc}
\begin{align*}
\int_{\lambda t}^t\frac{\phi'(s)}{\phi(s)}\hbox{d}s\geq\alpha \int_{\lambda t}^t\frac{\hbox{d}s}s,\\[-1.5pc]
\end{align*}
whence
\begin{align*}
\ln\frac{\phi(t)}{\phi(\lambda t)}\geq\ln\frac{t^{\alpha}}{(\lambda t)^{\alpha}}\\[-1.3pc]
\end{align*}
and consequently\vspace{-.2pc}
\begin{align*}
\phi(\lambda t)\leq\lambda^{\alpha}\phi(t).\\[-1.3pc]
\end{align*}
Let us take $t=t_0$. Then, for all $0<\lambda<1$, we have $\phi(\lambda
t_0)\leq\phi(t_0)\lambda^{\alpha}$, so $\phi(\lambda
t_0)\leq\frac{\phi(t_0)}{t_0^{\alpha}}\cdot (\lambda t_0)^{\alpha}$. If
we take $\epsilon=\alpha-1$, $A=\frac{\phi(t_0)}{t_0^{\alpha}}$ and
$u_0=t_0$, we get (b).\vspace{.4pc}

\noindent (b) $\Rightarrow$ (c). Take $\epsilon>0$, $A>0$ and $u_0>0$ such that
for all $0\leq u\leq u_0$, we have $\phi(u)\leq Au^{1+\epsilon}$. Since
$\frac 1 n\to0$ there exists $n_1\in\Bbb N$ such that $\frac 1 n\leq
u_0$ for all $n\geq n_1$. Therefore,\vspace{-.3pc}
\begin{align*}
\sum_{n=n_1}^{\infty}\phi\left(\frac{1}{n}\right)\leq \sum_{n=n_1}^{\infty}A\left(\frac{1}{n}\right)^{1+\epsilon}\leq
A\sum_{n=1}^{\infty}\frac 1 {n^{1+\epsilon}}<\infty.\\[-3.5pc]
\end{align*}
$\hfill\Box$\vspace{1pc}
\end{proof}

\begin{lemm}
{\it If $x\in l^0${\rm ,} $\{x_n\}\subset ces_{\phi}${\rm ,} $\sup\|x_n\|<\infty$ and \hbox{$0\leq x_n\!\uparrow\!x$} coordinatewise{\rm ,} then $x\in ces_{\phi}$ and $\|x_n\|\to\|x\|${\rm .}}
\end{lemm}

\begin{proof}
\looseness-1 Assume that $x_n\in ces_{\phi}$ for all $n \in \Bbb
N$, $\sup\|x_n\|<\infty$ and \hbox{$0\leq x_n(i)\!\uparrow\!x(i)$}
for each $i\in\Bbb N$. Denote $A=\sup_n\|x_n\|$. We know that
$\|x_n\|\leq A<\infty$ for all $n\in\Bbb N$, so $0\leq \frac
{x_n}A \leq\frac{x_n}{\|x_n\|}$ for all $n\in\Bbb N$. Therefore
$\rofi\left(\frac{x_n}{A}\right)\leq 1$ and since the modular
$\rofi$ is monotone,\break we get
\begin{equation*}
\rho_{\phi}\left(\frac{x_n}{A}\right)\leq \rho_{\phi}\left(\frac{x_n}{\|x_n\|}\right)\leq 1.
\end{equation*}
Then, by the Beppo Levi theorem and the fact that $A^{-1} x_n(i)\to
A^{-1} x(i)$ for each $i\in\Bbb N$, we get
\begin{equation*}
\rho_{\phi}\left(\frac{x}{A}\right)=\lim_{n\to\infty}\rho_{\phi}\left(\frac{x_n}{A}\right)=\sup_n\rofi\left(\frac{x_n}{A}\right)\leq 1,
\end{equation*}
\pagebreak

\noindent whence $x\in ces_{\phi}$ and $\|x\|\leq A$. By the
assumption that \hbox{$x_n\!\uparrow\!x$} coordinatewise and by
monotonicity of the norm, we get $\sup_n\|x_n\|\leq\|x\|$.
Therefore, we have $\|x\|=\sup_n\|x_n\|=\lim_{n\to\infty}\|x_n\|$.
$\hfill\Box$
\end{proof}

It is known that for any K\"{o}the sequence (function) space the Fatou
property implies its completeness (see \cite{17}). Therefore, $\cfi$ is a Banach space.

\begin{theor}[\!]
Let $A_{\phi}=\{x\in ces_{\phi}\!\!:\forall k>0\,\exists
n_k\,\sum_{n=n_k}^{\infty}\phi\left(\frac{k}{n}\sum_{i=1}^n|x(i)|\right)<\infty\}$. Then the following assertions are
true{\rm :}

\begin{enumerate}
\leftskip .3pc
\renewcommand{\labelenumi}{\rm (\roman{enumi})}
\item $A_{\phi}$ is a closed separable subspace of $\cfi${\rm ,}
\item $A_{\phi}=cl\{ x\in ces_{\phi}\!\!: x(i)\not=0 \text{ only for finite number of }\; i\in\Bbb N\}${\rm ,}
\item $A_{\phi}$ is the subspace of all order continous elements of $\cfi${\rm .}
\end{enumerate}
\end{theor}

\begin{proof}
It is easy to see that $A_\phi$ is a subspace of $\cfi$. Next we will
prove that $A_\phi$ is closed in $\cfi$. We must show that if $x_m\in
A_{\phi}$ for each $m\in\Bbb N$ and $x_m\to x \in ces_{\phi}$, then
$x\in A_{\phi}$. Take any $k>0$. We will show that there exists $n_k\in\Bbb N$ such that $\quad\sum_{n=n_k}^{\infty}\phi\left(\frac{k}{n}\sum_{i=1}^n|x(i)|\right)<\infty$. Since $\rho_{\phi}(k(x-x_m))\to 0$
for all $k>0$, there exists $M\in\Bbb N$ such that
$\rho_{\phi}(2k(x-x_M))<1$. Since $x_M\in A_{\phi}$, there exists
$n_M$ such that $\sum_{n=n_M}^{\infty}\phi
\left(\frac{2k}n\sum_{i=1}^n|x_M(i)|\right)<\infty$. As we will
see, we can take $n_k=n_M$. Indeed,
\begin{align*}
&\sum_{n=n_M}^{\infty}\phi\left(\frac k n\sum_{i=1}^n|x(i)|\right) = \sum_{n=n_M}^{\infty}\phi\left(\frac k
n\sum_{i=1}^n\left|\frac{2(x(i)-x_M(i))}2+\frac{2x_M(i)}2\right|\right)\\[.7pc]
&\quad\, \leq\sum_{n=n_M}^{\infty}\phi\left(\frac k n\sum_{i=1}^n\left|\frac{2(x(i)-x_M(i))}2\right|+\left|\frac{2x_M(i)}2\right|\right)\\[.7pc]
&\quad\,=\sum_{n=n_M}^{\infty}\phi\left(\frac 1 2\frac k n\sum_{i=1}^n |2(x(i)-x_M(i))|+\frac 1 2\frac k n\sum_{i=1}^n |2x_M(i)|\right)\\[.7pc]
&\quad\,\leq\sum_{n=n_M}^{\infty}\left(\frac 1 2\phi\left(\frac{2k}n\sum_{i=1}^n|x(i)-x_M(i)|\right)+\frac 1 2\phi\left(\frac{2k}n\sum_{i=1}^n|x_M(i)|\right)\right)\\[.7pc]
&\quad\,=\frac 1 2\sum_{n=n_M}^{\infty}\phi\left(\frac{2k}n\sum_{i=1}^n|x(i)-x_M(i)|\right)+ \frac 1 2\sum_{n=n_M}^{\infty}\phi\left(\frac{2k}n\sum_{i=1}^n|x_M(i)|\right)\\[.7pc]
&\quad\,\leq\frac 1 2\rofi(2k(x-x_M)+\frac 1
2\sum_{n=n_M}^{\infty}\phi\left(\frac{2k}n\sum_{i=1}^n|x_M(i)|\right)<\infty.
\end{align*}
By the arbitrariness of $k>0$, we get that $x\in A_{\phi}$, which
proves that $A_{\phi}$ is the closed subspace in the norm topology
in $\cfi$.

Now, we will prove assertion (ii). Let us define the set
$B_{\phi}=cl\{ x\in ces_{\phi}\!\!: x(i)=0 \text{ for a.e. $i\in\Bbb
N$}\}$. We will prove that $A_{\phi}$ and $B_{\phi}$ are equal.

First we will show that $B_{\phi}\subset A_{\phi}$. If
$B_{\phi}=\emptyset$, the inclusion $B_{\phi}\subset A_{\phi}$ is
obvious. So, assume that $B_{\phi}\not=\emptyset$. Take
$x=(\underbrace{0,\dots,0}_{l-1\,{\rm times}},1,0,0,\dots)\in
B_{\phi}$ and $k>0$. We have from Theorem~2.1 that there exists
$n_k$ such that\pagebreak
\begin{equation*}
\sum\limits_{n=n_k}^{\infty}\phi\left(\frac k n\right)<\infty.
\end{equation*}
We can assume that $n_k\geq l$. Hence $x\in A_{\phi}$, and
so, by the fact that $A_{\phi}$ is a linear subspace of $\cfi$, we
get the inclusion $B_{\phi}\subset A_{\phi}$.

Now, we will show that $A_{\phi}\subset B_{\phi}$. Let
$x=(x_1,x_2,\dots,x_k,x_{k+1},\dots)\in A_{\phi}$ and define
$x^k=(x_1,x_2,\dots,x_k,0,0,\dots)$ for any $k\in\Bbb N$. Obviously
$x^k\in B_{\phi}$. We will show that
$\rho_{\phi}(\alpha(x-x_k))\rightarrow 0$ for each $\alpha>0$.
Take any $\alpha>0$ and $\epsilon>0$. Since $x\in A_{\phi}$, so
there exists $k_0\in\Bbb N$ such that
\begin{equation*}
\sum\limits_{n=k_0+1}^{\infty}\phi\left(\frac {\alpha} n\sum\limits_{i=1}^n|x(i)|\right)<\epsilon.
\end{equation*}
Then for any $k\geq k_0$,
\begin{align*}
\rho_{\phi}(\alpha(x-x^k))&\leq\rho_{\phi}(\alpha(x-x^{k_0}))=\rho_{\phi}(\alpha(0,\dots,0,x_{k_0+1},x_{k_0+2},\dots))\\[.3pc]
&=\sum\limits_{n=k_0+1}^{\infty}\phi\left(\frac {\alpha}
n\sum\limits_{i=k_0+1}^n|x(i)|\right)\\[.6pc]
&\leq \sum\limits_{n=k_0+1}^{\infty}\phi\left(\frac {\alpha}
n\sum\limits_{i=1}^n|x(i)|\right)<\epsilon.
\end{align*}

Next we will prove assertion (iii). Let $x\in A_{\phi}$. We will
show that $x$ is order continuous.  Take any $k>0$ and
$\epsilon>0$. Then there exists $n_k\in\Bbb N$ such that
$\sum_{n=n_k}^{\infty}\phi\left(\frac k
n\sum_{i=1}^n|x(i)|\right)<\frac{\epsilon}2$. Assume that
$x_m\downarrow 0$ coordinatewise and $x_m\leq |x|$ for all $m \in \Bbb N$. Denote
\begin{equation*}
\phi\left(\frac k n\sum_{i=1}^n|x(i)|\right)=\alpha(n)
\end{equation*}
and
\begin{equation*}
\phi\left(\frac k n\sum_{i=1}^n|x_m(i)|\right)=\alpha_m(n)
\quad\hbox{for any}\,n\in\Bbb N.
\end{equation*}
Since $x_{m}\downarrow 0$ coordinatewise, we get
$\alpha_m(n)\rightarrow 0$ as $m\rightarrow \infty$ for any
$n\in\Bbb N$. Consequently, there is $m_{\epsilon}\in\Bbb N$ such
that $\sum_{n=1}^{n_k-1}\alpha_m(n)<\frac{\epsilon}2$ for any
$m\geq m_{\epsilon}$. Moreover,
$\sum_{n=n_k}^{\infty}\alpha_m(n)<\sum_{n=n_k}^{\infty}\alpha(n)<\frac{\epsilon}2$
for all $n\geq n_k$ and $m\in \Bbb N$. Therefore
$\rofi(kx_m)<\epsilon$ for all $m\geq m_{\epsilon}$, which means
that $\rofi(kx_m)\rightarrow 0$. By the arbitrariness of $k>0$,
this means that\break $\|x_m\|\rightarrow 0$.

Let $x\in\cfi$ be an order continuous element. Since
\begin{equation*}
\|(0,\dots,0,x(n+1),x(n+2),\dots)\|\rightarrow
0\quad\text{as}\, n\rightarrow \infty,
\end{equation*}
so it easy to see that $x\in cl\{ x\in ces_{\phi}\!\!: x(i)=0 \text{
for a.e. $i\in\Bbb N$}\}$.

Finally, we will show that $A_\phi$ is separable. Roughly
speaking, this follows by the fact that the counting measure on
$\Bbb N$ is separable and $A_\phi$ is order continuous.

Define the set $C_{\phi}=cl\{ x\in ces_{\phi}\!\!: x(i)=0$ for a.e.
$i\in\Bbb N$  and $x(i)\in Q\}$ which is countable. It is obvious
that $C_{\phi}\subset B_{\phi}$. Now, we will show that
$B_{\phi}\subset C_{\phi}$. Let
$x=(x(1),x(2),\dots,x(k),0,0,\dots)\in B_{\phi}$ and $
x_m=(x_m(1),\dots,x_m(k),0,\dots)\in C_{\phi}$ will be such that
$x_m(i)\rightarrow x(i)$ as $m\rightarrow \infty$. We will show
that $\|x_m-x\|\rightarrow 0$.

Let us take any $\lambda>0$. We have
\begin{equation*}
\lambda (|x(1)-x_m(1)| + |x(2)-x_m(2)| + \cdots +|x(k)-x_m(k)|)
\leq 1
\end{equation*}
for $m$ large enough. Then by convexity of $\phi$,
\begin{align*}
\hskip -4pc \rofi(\lambda(x-x_m)) &\leq \sum_{n=1}^{\infty} \phi
\left( \lambda \frac{ |x(1)-x_m(1)| + |x(2)-x_m(2)| + \cdots
+|x(k)-x_m(k)|
} n \right)\\[.5pc]
&\leq \lambda (|x(1)-x_m(1)| + |x(2)-x_m(2)| + \cdots
+|x(k)-x_m(k)|)\\[.5pc]
&\quad\, \times \sum_{n=1}^{\infty} \phi \left( \frac 1 n
\right)\rightarrow 0
\end{align*}
as $m\rightarrow \infty$. By the arbitrariness of $\lambda$, we
have $\|x_m-x\|\rightarrow 0$ as $m\rightarrow \infty$.
Consequently, $B_{\phi}=C_{\phi}$. Since $B_{\phi}=A_{\phi}$ and
the space $C_{\phi}$ is separable, we get the separability of
$A_{\phi}$.

\hfill $\Box$\vspace{.3pc}
\end{proof}

\begin{theor}[\!]
If $\phi\in\Delta_2(0)${\rm ,} then $A_{\phi}=ces_{\phi}$.\vspace{.2pc}
\end{theor}

\begin{proof}
We should only show that $ces_{\phi}\subset A_{\phi}$. Let $x\in
ces_{\phi}$. Then there exists $\alpha>0$ such that
$\rho_{\phi}(\alpha x)<\infty$. We will show that for any
$\lambda>0$ there exists $n_{\lambda}$ such that
$\sum_{n=n_{\lambda}}^{\infty}\phi\left(\frac{\lambda}n\sum_{i=1}^n|x(i)|\right)<\infty$.
We take only $\lambda>\alpha$, because for $\lambda<\alpha$ we
have $\sum_{n=n_1}^{\infty}\phi\left(\frac{\lambda}n
\sum_{i=1}^n|x(i)|\right)<\sum_{n=n_1}^{\infty}\phi
\left(\frac{\alpha}n\sum_{i=1}^n|x(i)|\right)<\infty$ from
monotonicity of the function $\phi$. Let $\lambda>\alpha$. By
$\phi\in\Delta_2(0)$, we have that $\phi\in\Delta_l(0)$ for any
$l>1$, whence for $l:=\frac{\lambda}{\alpha}$ there exists $k,
u_0>0$ such that $\phi(lu)\leq k\phi(u)$ for all $u\leq u_0$. By
$\rofi(\alpha x)<\infty$, there exists $n_{\lambda}$ such that
$\frac{\alpha}n\sum_{i=1}^n|x(i)|<u_0$ for all $n\geq
n_{\lambda}$. Therefore,
\begin{align*}
\sum_{n=n_{\lambda}}^{\infty}\phi\left(\frac{\lambda}n\sum_{i=1}^n|x(i)|\right)
& = \sum_{n=n_{\lambda}}^{\infty}\phi\left(
\frac{\lambda\alpha}{\alpha n}\sum_{i=1}^n|x(i)|\right)\\[.5pc]
&\leq
k\sum_{n=n_{\lambda}}^{\infty}\phi\left(\frac{\alpha}n\sum_{i=1}^n|x(i)|\right)<\infty,
\end{align*}
and the proof is finished.\hfill $\Box$
\end{proof}

\setcounter{theore}{0}
\begin{coro}$\left.\right.$\vspace{.5pc}

\noindent If $\phi\in\Delta_2(0)${\rm ,} then
\begin{enumerate}
\renewcommand\labelenumi{{\rm (\roman{enumi})}}
\leftskip .1pc
\item the space $\cfi$ is a separable{\rm ,}
\item the space $\cfi$ is order continuous.
\end{enumerate}
\end{coro}

We will assume in the following that the function $\phi$ is
finite. We will prove some useful lemmas.

\setcounter{theore}{0}
\begin{lem}
For any $x\in A_{\phi}${\rm ,}
\begin{equation*}
\|x\|=1\quad\hbox{if and only if}\,\rofi(x)=1.
\end{equation*}
\end{lem}
\pagebreak

\begin{proof}
We need only to show that $\|x\|=1$ implies $\rofi(x)=1$ because
the opposite implication holds in any modular space. Assume that
$\phi<\infty$ and take $x\in A_{\phi}$ with $\|x\|=1$. Note that
$\rofi(x)\leq1$. Assume that $\rofi(x)<1$. Since $x\in A_{\phi}$,
we have that  $\rofi(kx)<\infty$ for all $k>0$. Let us define the
function $f(\lambda)=\rofi(\lambda x)$, which is convex and  has
finite values. Hence $f$ is continous on $\Bbb R_{+}$ and $f(1)<1$
by the assumption that $\rofi(x)<1$. Then, by the continuity of
$f$ there exists $r>1$ such that $f(r)\leq 1$, that is,
$\rho_{\phi}(rx)\leq 1$. Then $\|rx\|\leq 1$, whence
$\|x\|\leq\frac 1 r<1$, a contradiction, which shows that
$\rofi(x)=f(1)=1$. \hfill $\Box$
\end{proof}

\begin{lem}
If $\phi\in\Delta_2(0)${\rm ,} then $\rofi\in \Delta_2^S$.
\end{lem}

\begin{proof}
Take arbitrary $\epsilon > 0$, $a>0$ and $\rofi(x)\leq a$. Then
$\rofi(x)=\sum_{n=1}^{\infty}\phi\left(\sigma x(n)\right)\leq a$,
whence $\phi(\sigma x(n))\leq a$ for any $n\in\Bbb N$. If $b>0$ is
the number satisfying $\phi(b)=a$, then $\sigma x(n)\leq b$ for
any $n\in\Bbb N$. Since $\phi\in\Delta_2(0)$ and $\phi<\infty$, so
$\phi\in\Delta_2([0,b])$, i.e. there exists $K>0$ such that
$\phi(2u)\leq K\phi(u)$ for all $u\in[0,b]$. We have
\begin{align*}
\rofi(2x)&=\sum_{n=1}^{\infty}\phi(\sigma 2x(n))=\sum_{n=1}^{\infty}\phi(2\sigma x(n))\\[.5pc]
&\leq k\sum_{n=1}^{\infty}\phi(\sigma x(n))=k\rofi(x).
\end{align*}
$\left.\right.$\vspace{-3pc} \hfill $\Box$
\end{proof}\vspace{1.5pc}

\begin{lem}
Assume that $\phi\in\Delta_2(0)$. Then for any $L>0$ and $\epsilon
>0$ there exists $\delta=\delta(L,\epsilon)>0$ such that
\begin{equation*}
|\rho_{\phi}(x+y)-\rho_{\phi}(x)|<\epsilon
\end{equation*}
for all $x,y \in ces_{\phi}$ with $\rho_{\phi}(x)\leq L$ and
$\rho_{\phi}(y)\leq \delta(L,\epsilon)$.
\end{lem}

\begin{proof}
In virtue of Lemma 2.2 it suffices to apply Lemma 2.1 in \cite{4}.
\hfill $\Box$
\end{proof}


\begin{lem}
If $\phi\!\in\!\Delta_2(0)${\rm ,} then for any sequence $(x_n)\in ces_{\phi}$
the condition$\|x_n\| \rightarrow 0$ holds if and only if $\rho_{\phi}(x_n)
\rightarrow 0$.
\end{lem}

\begin{proof}
It suffices to apply Lemmas~2.2 and 2.3 in \cite{4}. \hfill $\Box$
\end{proof}

\begin{lem}
If $\phi\in\Delta_2(0)${\rm ,} then for any $x\in ces_{\phi}${\rm ,}
\begin{equation*}
\|x\| = 1\,\textit{if and only if}\,\rho_{\phi}(x) = 1.
\end{equation*}
\end{lem}

\begin{proof}
The result follows from Lemma~2.2 and Corollary~2.2 in \cite{4}.
\hfill $\Box$
\end{proof}

\begin{lem}
If $\phi\in\Delta_2(0)${\rm ,} then for any $\epsilon > 0$ there exists
$\delta = \delta (\epsilon) > 0$ such that $\|x\| \geq 1 + \delta$ whenever
$x \in \cfi$ and $\rho_{\phi}(x) \geq 1 +\epsilon$.
\end{lem}

\begin{proof}
The result follows by applying Lemmas~2.2 and 2.4 in \cite{4}.
\hfill $\Box$
\end{proof}

\begin{lem}
Let $\phi \in \Delta_2 (0)$. Then for each $\epsilon > 0$ there exists
$\delta = \delta (\epsilon)$ such that $\rofi (x) > \delta$ whenever
$\|x\| \geq \epsilon$.
\end{lem}

\begin{proof}
Suppose for the contrary there exists $\epsilon>0$ such that for any
$\delta > 0$, there exists $x$ such that $\rho_{\phi} (x) \leq \delta$ and
$\|x\| \geq \epsilon$. Take $\delta_n = \frac 1 n$ and the sequence
$(x_n)_{n \in \Bbb N}$ in $\cfi$ satisfying $\rofi (x_n) \leq \frac 1 n$ and
$\|x_n\| \geq \epsilon$. Consequently $\rofi(x_n) \rightarrow 0$ as $n \rightarrow \infty$.
From Lemma~2.4 it follows that $\|x_n\| \rightarrow 0$, a contradiction finishing
the proof.
\hfill $\Box$
\end{proof}

\begin{lem}
If $\phi\in\Delta_2(0)${\rm ,} then $\|x_n\| \rightarrow \infty$ whenever
$\rho_{\phi}(x_n) \rightarrow \infty$.
\end{lem}

\begin{proof}
Suppose $(\|x_n\|)$ is a bounded sequence, that is, there exists $M > 0$
such that $\|x_n\|\leq M$ for all $n\in\Bbb N$. Take $s\in\Bbb N$ such
that $M\leq 2^s$. Then $\|x_n\|\leq 2^s$, whence
$\|\frac{x_n}{2^s}\|\leq 1$ and
$\rho_{\phi}\left(\frac{x_n}{2^s}\right)\leq 1$. Consequently,
$\phi((\sigma\frac{x_n}{2^s})(i))\leq 1$ for all
$i\in\Bbb N$, and then, there exists some $L>0$ such that
$\left(\sigma\frac{x_n}{2^s}\right)(i)\leq L$ for all $i\in\Bbb N$.
Since $\phi\in\Delta_2(0)$ and $\phi<\infty$,
$\phi\in\Delta_2([0,2^{s-1} L])$. We have for all $n\in\Bbb N$,
\begin{equation*}
\rho_{\phi}(x_n)=\rho_{\phi}\left(2^s\frac{x_n}{2^s}\right)
\leq k^s\rho_{\phi}\left(\frac{x_n}{2^s}\right)\leq k^s,
\end{equation*}
whence $\rho_{\phi}(x_n)\not \rightarrow \infty$.
\hfill $\Box$\vspace{.3pc}
\end{proof}

\begin{lem}
If $\phi\in\Delta_2(0)${\rm ,} then for any sequence $(x_n)$ in
$\cfi${\rm ,} we have
\begin{equation*}
\|x_n\| \rightarrow 1 \quad \textit{if and only if}\,\rofi(x_n)
\rightarrow 1.
\end{equation*}
\end{lem}

\begin{proof}
The implication $\rofi(x_n) \rightarrow 1 \Rightarrow \|x_n\|
\rightarrow 1$ is almost obvious. Namely, we have $\rofi(x)\leq\|x\|$ if
$\rofi(x)\leq 1$ and $\|x\|\leq\rofi(x)$ if $\rofi(x)>1$. Therefore
$|\|x_n\|-1|\leq|\rofi(x_n)-1|$ and the result follows.
Now, assuming that $\|x_n\| \rightarrow 1$, we consider two cases:

\begin{enumerate}
\renewcommand\labelenumi{\arabic{enumi}.}
\leftskip .15pc
\item $\|x_n\|\uparrow 1$. From Lemma~2.8 we know that the sequence
$(\rofi(2x_n))$ is bounded, that is, there exists $A>0$ such that
$\rofi(2x_n)\leq A$ for all $n\in\Bbb N$. Assume for the contrary that
$\rofi(x_n)\not \rightarrow 1$. We can assume that $\|x_n\|>\frac 1 2$
for all $n\in\Bbb N$ and there exists $\epsilon>0$ such that
$\rofi(x_n)<1-\epsilon$ for all $n\in\Bbb N$. Take $a_n:=\frac
1{\|x_n\|}-1$. Then $a_n \rightarrow 0$ and $a_n\leq 1$. By Lemma 2.5, we
have
\begin{align*}
1&=\rofi\left(\frac{x_n}{\|x_n\|}\right)=\rofi((a_n+1)x_n)\\[.3pc]
&=\rofi(2a_nx_n+(1-a_n)x_n)\leq a_n\rofi(2x_n)+(1-a_n)\rofi(x_n)\\[.3pc]
&\leq a_n\cdot A+(1-a_n)(1-\epsilon) \rightarrow 1-\epsilon
\end{align*}
as $n \rightarrow \infty$, a contradiction.

\item $\|x_n\|\!\downarrow\!1$. Assume that $\|x_n\|\leq 2$ for $n\in\Bbb
N$ and there exists $\epsilon>0$ such that $\rofi(x_n)>1+\epsilon$ for
all $n\in\Bbb N$. From Lemma 2.8 we know that there exists $B>0$ such
that $\rofi(2x_n)\leq B$ for all $n\in\Bbb N$. By the assumption we have
$0\leq 1-\frac 1{\|x_n\|}\leq 1$, $0\leq 2-\|x_n\|\leq 1$. The
inequality $\frac 1 a +a\geq 2$ for any $a>0$ yields $0\leq \big(1-\frac
1{\|x_n\|} \big)+(2-\|x_n\|)=3- \big(\frac 1{\|x_n\|}+\|x_n\|\big)\leq 3-2=1$ for any
$n\in\Bbb N$. Therefore, we have
\begin{align*}
&1+\epsilon\leq \rofi(x_n)=\rofi\left(\left(1-\frac 1{\|x_n\|}\right)\cdot 2x_n+(2-\|x_n\|)\frac {x_n}{\|x_n\|}\right)\\[.3pc]
&\qquad\,\leq \left(1-\frac 1{\|x_n\|}\right)\rofi(2x_n)+(2-\|x_n\|)\rofi\left(\frac {x_n}{\|x_n\|}\right)\\[.3pc]
&\qquad\,\leq\left(1-\frac 1{\|x_n\|}\right)B+\rofi\left(\frac {x_n}{\|x_n\|}\right) \rightarrow 1,
\end{align*}
because $\rofi \big(\frac {x_n}{\|x_n\|}\big)=1$ for any $n\in\Bbb N$ and
$1-\frac 1 {\|x_n\|} \rightarrow 0$, a contradiction which finishes the
proof.\hfill$\Box$\vspace{-1.5pc}
\end{enumerate}
\end{proof}

Now we will consider monotonicity properties of $A_\phi$ and $\cfi$.

\setcounter{theore}{4}
\begin{theor}[\!]
The space $A_{\phi}$ is strictly monotone if and only if $\phi>0$.
\end{theor}

\begin{proof}
Denote $a_\phi=\sup\set{t\geq 0\!\!: \phi(t)=0}$ and assume that $a_\phi>0$.
We will show that under this assumption there exists $x,y\in ces_{\phi}$
such that $x\leq y$, $x\not=y$ and $\|x\|=\|y\|$. We define the function
$f(t)=\sum_{n=1}^{\infty}\phi\left(\frac t n\right)$ for $t\geq 0$.
Since $a_{\phi}>0$, $\frac t n \rightarrow 0$ as $n \rightarrow \infty$
and $a_\phi>0$, so $\sum_{n=1}^{\infty}\phi\left(\frac t n\right)$ is
convergent for all $t\in \Bbb R_+$. Since $\phi$ is a convex function,
so $ f$ is convex, too. Then $f$ is continuous on $\Bbb R_+$ and $f(t)
\rightarrow \infty$ as $t \rightarrow \infty$, whence $f(\Bbb R_+)=\Bbb
R_+$ and by the Darboux property of $f$ we know that there exists $c\in
\Bbb R$ such that $f(c)=\sum_{n=1}^{\infty}\phi\left(\frac c
n\right)=1$. Since $\frac{c+1} n \rightarrow 0$ as $n \rightarrow
\infty$, there exists $n_0$ such that $\frac {c+1}{n_0}\leq
a_{\phi}$. Consider two sequences $x=(c,0,0, \dots)$ and
$y=(\underbrace{c,0, \dots,0}_{n_0-1\text{ times}},1,0, \dots)$. It is obvious
that $x\not= y$ and $x<y$. Moreover,
\begin{align*}
\rho_{\phi}(x) &=\phi(c)+\phi\left(\frac c 2\right)+\phi\left(\frac c
3\right)+ \cdots =f(c)=1,\\[.5pc]
\rho_{\phi}(y) &=\phi(c)+\phi\left(\frac c 2\right)+ \cdots
+\phi\left(\frac{c}{n_0-
1}\right)+\phi\left(\frac{c+1}{n_0}\right)\\[.5pc]
&\quad\, +\phi\left(\frac{c+1}{n_0+1} \right) + \cdots = 1.
\end{align*}
Since $\rho_{\phi}(x)=\rho_{\phi}(y)=1$, we have $\|x\|=\|y\|=1$, which means
that $A_\phi\notin (SM)$.

Assume now that $a_\phi=0$, $y\geq x\geq 0$, $x\not=y$ and $x,y\in
A_{\phi}$. We can assume that $\|x\|=1$. From Lemma 2.1 we know that
$\rho_{\phi}(x)=1$. In order to show that $\|y\|>1$ we need to show that
$\rho_{\phi}(y)>1$. Note that $\rofi(x+y)\geq \rofi(x)+\rofi(y)$ for all
nonnegative $x,y\in A_{\phi}$. Therefore
\begin{equation*}
\rho_{\phi}(y)=\rho_{\phi}(x+(y-x))\geq\rho_{\phi}(x)+\rho_{\phi}(y-x)=
1+\rho_{\phi}(y-x)>1,
\end{equation*}
because of $y-x>0$ and $\phi>0$, whence $\rofi(y-x)>0$. This finishes the proof.
\hfill $\Box$
\end{proof}

From the last theorem, we get the following.

\setcounter{theore}{1}
\begin{coro}$\left.\right.$\vspace{.5pc}

\noindent If the space $\cfi$ is strictly monotone{\rm ,} then $\phi>0$.
\end{coro}

Before formulating the next theorem note that $\phi>0$ whenever
$\phi\in\Delta_2(0)$.

\setcounter{theore}{5}
\begin{theor}[\!]
If $\phi\in\Delta_2(0)${\rm ,} then $\cfi$ is uniformly monotone.
\end{theor}

\begin{proof}
Let $\epsilon>0$ and $x,y\geq 0$ be such that $\|x\|=1$ and
$\|y\|\geq\epsilon$. From Lemma 2.5 we have $\rofi(x)=1$ and from Lemma
2.7 we have that $\rofi(y)>\eta$ where $\eta>0$ is independent of $y$.
Then
\begin{equation*}
\rofi(x+y)\geq\rofi(x)+\rofi(y)\geq1+\eta.
\end{equation*}
By Lemma~2.6, there exists $\delta > 0$ independent of $x$ and $y$ such that  $\|x+y\|\geq1+\delta$.

\hfill $\Box$
\end{proof}

Next we consider rotundity of $\cfi$. In order to be able to prove
criteria for rotundity of $\cfi$, we need first to prove the following.

\setcounter{theore}{9}
\begin{lem}
Let $\phi\in\Delta_2(0)$ and $y, z\in S(\cfi)$ satisfy $\frac
{y+z}{2}\in S(\cfi).$ If $y\not= z${\rm ,} then there exists $i_0\in\Bbb N$
such that $|y(i_0)|\not= |z(i_0)|$.
\end{lem}

\begin{proof}
Assume for the contrary that the assumptions are satisfied, $y\not=z$ and $|y|=|z|$. Then there is
$i_0\!\in\!\Bbb N$ such that $y(i_0)\!\not\!=z(i_0)$, but $|y(i_0)|=|z(i_0)|$,
whence \hbox{$y(i_0)+z(i_0)=0$.} Consequently,
\begin{align*}
1&=\rofi\left(\frac{y+z}2\right)=\sum_{n=1}^{\infty}\phi\left(\frac 1 n\sum_{i=1}^n\frac{|y(i)+z(i)|}2\right)\\[.3pc]
&=\sum_{n=1}^{\infty}\phi\left(\frac 1 2\sum_{i=1}^n\frac{|y(i)+z(i)|}n\right)=\sum_{n=1}^{\infty}\phi\left(\frac 1 2\sum_{i\in\Bbb N\setminus \{i_0\}}\frac{|y(i)+z(i)|}n\right)\\[.3pc]
&\leq\sum_{n=1}^{\infty}\phi\left(\frac 1 2\left(\frac 1 n\sum_{i\in N\backslash \{i_0\}}|y(i)|+\frac 1 n\sum_{i\in N\backslash \{i_0\}}|z(i)|\right)\right)\\[.3pc]
&\leq\sum_{n=1}^{\infty}\left(\frac 1 2\phi\left(\frac 1 n\sum_{i\in N\backslash \{i_0\}}|y(i)|\right)+\frac 1 2\phi\left(\frac 1 n\sum_{i\in N\backslash \{i_0\}}|z(i)|\right)\right)\\[.3pc]
&<\frac 1 2\rofi(y)+\frac 1 2\rofi(z)=1,
\end{align*}
a contradiction which finishes the proof.
\hfill $\Box$
\end{proof}

Given any Orlicz function $\phi$ with values in $\Bbb R_+$ such that
$\sum_{i=1}^\infty\phi\left(\frac 1 i\right)<\infty$, define the
function
\begin{equation}
f(a)=2\phi(a)+\sum_{i=3}^{\infty}\phi\left(\frac{2} i a\right).
\end{equation}
Since the function $\phi$ is convex, so $f$ is convex as well. By
Theorem~2.1 it has finite values. Therefore $f$ is continuous and $f(a)
\rightarrow \infty$ as $a \rightarrow \infty$, whence we deduce that
there exists $\alpha\in\Bbb R$ such that $f(\alpha)=1$.


\setcounter{theore}{6}
\begin{theor}[\!]
If $\phi\in\Delta_2(0)$ then $\cfi$ is rotund if and only if $\phi$ is
strictly convex on the interval $[0,\alpha]${\rm ,} where $f(\alpha)=1$ and
$f$ is defined by formula (2.1).
\end{theor}

\begin{proof}

Suppose $\phi$ is not strictly convex on $[0, \alpha ]$. Then there
exists an interval $[b, c]\subset (0, \alpha)$ on which $\phi$ is
affine.

Since $c<\alpha$, we have
\begin{equation*}
2\phi(c)+\sum_{i=3}^\infty \phi\left(\frac{2c}{i} \right)<1.
\end{equation*}
Take $d>0$ such that
\begin{equation*}
2\phi(c)+\sum_{i=3}^\infty \phi\left(\frac{2c+d}{i} \right)<1.
\end{equation*}
Choose $b_1, c_1$ such that $b<b_1<c_1<c$ and
\begin{align*}
\phi(b)+\phi\left(\frac{b+c}{2}\right) &=
\phi(b_1)+\phi\left(\frac{b_1+c_1}{2}\right),\\[.5pc]
&\quad\,b_1-b<\frac{d}{2}\quad \text {and} \quad c-c_1<\frac{d}{2}.
\end{align*}
By  $|b+c-b_1-c_1|<d,$ there is $k>0$ for which either $b+c=b_1+c_1+k$ or $b+c+k=b_1+c_1.$

Without loss of generality, we may assume that $b+c+k=b_1+c_1,$ whence
\begin{align*}
\phi(b) &+\phi \left(\frac{b+c}{2}\right)+\sum_{i=3}^\infty
\phi\left(\frac{b+c+k}{i}\right)\\[.7pc]
&\quad\, = \phi(b_1)+\phi\left(\frac{b_1+c_1}{2}\right)
 +\sum_{i= 3}^\infty \phi\left(\frac{b_1+c_1}{i}\right).
\end{align*}
Take $k_1>0$ such that
\begin{equation}
\phi(b)+\phi\left(\frac{b+c}{2}\right)+\phi\left(\frac{b+c+k}{3}\right)
+\sum_{i=4}^\infty \phi\left(\frac{b+c+k+k_1}{i}\right)=1.
\end{equation}
Since $b+c+k=b_1+c_1,$ we have
\begin{equation}
\phi(b_1)+\phi\left(\frac{b_1+c_1}{2}\right)+\phi\left(\frac{b_1+c_1}{3}
\right)+\sum_{i=4}^\infty \phi\left(\frac{b_1+c_1+k_1}{i}\right)=1.
\end{equation}
Put
\begin{equation*}
x=(b, c, k, k_1, 0, 0,\dots)
\end{equation*}
and
\begin{equation*}
y=(b_1, c_1, 0, k_1, 0, 0,\dots).
\end{equation*}
By (2.2) and (2.3), we have $\rofi(x)=1=\rofi(y)$. So, Lemma~2.5 yields
$x, y\in S(\cfi).$ Again, by (2.2) and (2.3) and the fact that
$\phi$ is affine on $[b, c]$, we have
\begin{align*}
\rofi\left(\frac{x+y}{2}\right) &=
\phi\left(\frac{b+b_1}{2}\right)+\phi\left(\frac{\frac{b+c}{2}+\frac{b_1+c_1}{2}}{2}\right)+\phi
\left(\frac{b+c+k}{3}\right)\\[.5pc]
&\quad +\sum_{i=4}^\infty \phi\left(\frac{b+c+k+k_1}{i}\right)\\[.5pc]
&= \frac{1}{2}\left(\phi(b)+\phi(b_1)\right)+\frac{1}{2}\left(\phi\left(\frac{b+c}{2}\right)+\phi\left(
\frac{b_1+c_1}{2}\right)\right)\\[.5pc]
&\quad + \phi\left(\frac{b+c+k}{3}\right) + \sum_{i=4}^\infty \phi\left(\frac{b+c+k+k_1}{i}\right) = 1.
\end{align*}
Therefore Lemma~2.5 yields $\big\|\frac{x+y}{2}\big\|=1$, which means that
$ces_\phi$ is not rotund.\pagebreak

Conversely, let $x\in S(ces_\phi)$. We need to prove that $x$ is an
extreme point. If $x$ is not an extreme point, then there exists $y, z\in
S(ces_\phi)$ such that $2x=y+z$ and $y\not= z.$ We will prove that
$|y|=|z|$ and by Lemma 2.10, we will get a contradiction, finishing the
proof.

Since $\phi\in\Delta_2(0),$ Lemma 2.5 yields that
$\rofi(x)=\rofi(y)=\rofi(z)=1$ and
\begin{align*}
1=\rofi(x)= \rofi\left(\frac {y+z}{2}\right) &= \sum_{n=1}^\infty
\phi\left(\frac{1}{n}\sum_{i=1}^n \frac {|y(i)+z(i)|}{2}\right)\\[.5pc]
&\leq\sum_{n=1}^\infty \phi\left(\frac{1}{n}\sum_{i=1}^n \frac
{|y(i)|+|z(i)|}{2}\right)\\[.5pc]
&\leq\frac{1}{2}\left[\sum_{n=1}^\infty
\phi\left(\frac{1}{n}\sum_{i=1}^n|y(i)|\right)\right.\\[.3pc]
&\hskip 1cm \left. +\sum_{n=1}^\infty
\phi\left(\frac{1}{n}\sum_{i=1}^n|z(i)|\right)\right]\\[.5pc]
&= \frac {1}{2}[\rofi(y)+\rofi(z)]\\[.3pc]
&= 1.
\end{align*}
Thus for each $n\in\Bbb N$ we have
\begin{equation}
\phi\left(\frac{1}{n}\sum_{i=1}^n \frac
{|y(i)|+|z(i)|}{2}\right)=\frac{1}{2}\left[\phi\left(\frac{1}{n}\sum_{i=
1}^n|y(i)|\right)+\phi\left(\frac{1}{n}\sum_{i=1}^n|z(i)|\right)\right].
\end{equation}

\begin{case}
{\rm $\frac{1}{n}\sum_{i=1}^n|x(i)|\leq \alpha$ for each $n\in \Bbb N$. By condition (2.4) and the fact that  $\phi$ is strictly convex on the interval $[0, \alpha]$, we have
$\frac 1 n\sum_{i=1}^n |y(i)|=\frac 1 n\sum_{i=1}^n|z(i)|$ for each $n\in \Bbb N$. Consequently{\rm ,} $|y|=|z|$.}
\end{case}

\begin{case}
{\rm There exists $n$ such that $\frac{1}{n}\sum_{i=1}^{n}|x(i)|> \alpha$.
We claim that there exists only one such $n$. Assume for the contrary
that there exists $n_0<n_1$ such that
$\frac{1}{n_0}\sum_{i=1}^{n_0}|x(i)|> \alpha$ and
$\frac{1}{n_1}\sum_{i=1}^{n_1}|x(i)|> \alpha$. Then $n_1\geq 2$ and we
have
\begin{align*}
1 &= \rofi(x)>2\phi(\alpha)+\sum_{i=n_1+1}^\infty \phi\left(\frac {n_1\alpha}i\right)=
2\phi(\alpha)+\sum_{i=1}^\infty \phi\left(\frac {n_1\alpha}{n_1+i}\right)\\[.7pc]
&\geq 2\phi(\alpha)+\sum_{i=1}^\infty \phi\left(\frac {2\alpha}{2+i}\right)=
2\phi(\alpha)+\sum_{i=3}^\infty \phi\left(\frac {2\alpha}{i}\right)=1,
\end{align*}
a contradiction, which proves the Claim. Let $n_0$ be the only natural
number for which $\frac{1}{n_0}\sum_{i=1}^{n_0}|x(i)|> \alpha$.
As in Case I{\rm ,} we can prove that
$\frac 1 n\sum_{i=1}^n |y(i)|=\frac 1 n\sum_{i=1}^n|z(i)|$ for each $n\not=n_0$. Since $\rofi(y)=\rofi(z)=1${\rm ,} we get\pagebreak
\begin{align*}
\phi\left(\frac 1 {n_0}\sum_{i=1}^{n_0} |y(i)|\right) &=
1-\sum_{n\in\Bbb N\setminus\set{n_0}}\phi\left(\frac 1 n \sum_{i=1}^n |y(i)|\right)\\[.7pc]
& = 1-\sum_{n\in\Bbb N\setminus\set{n_0}}\phi\left(\frac 1 n \sum_{i=1}^n |z(i)|\right)=
\phi\left(\frac 1 {n_0}\sum_{i=1}^{n_0} |z(i)|\right).
\end{align*}
Consequently, $|y|=|z|$. This finishes the proof. \hfill$\Box$}
\end{case}\vspace{-1pc}
\end{proof}

\setcounter{theore}{0}
\begin{rem}
{\rm Note that criteria for rotundity of Ces\`{a}ro--Orlicz sequence spaces
$\cfi$ are weaker than criteria for rotundity of Orlicz sequence spaces
$l_{\phi}$. Namely, we can easily conclude from \cite{11} that an Orlicz
sequence space $l_{\phi}$ is rotund if and only if $\phi$ attains value
1, $\phi\in\Delta_2(0)$ and $\phi$ is strictly convex on the interval
$[0,a]$ where $\phi(a)=\frac 1 2$, which is smaller from the interval
$[0,\alpha]$, where $\alpha$ is defined by (2.1).}
\end{rem}

\end{document}